\documentclass[11pt]{article}

\title{A counterexample to a conjecture of S.E. Morris}
\author{J. F. Feinstein}

\usepackage{amsfonts}
\usepackage{latexsym}
\usepackage{theorem}

\theorembodyfont{\sl}
\newtheorem{theorem}{Theorem}[section]
\newtheorem{proposition}[theorem]{Proposition}
\newtheorem{lemma}[theorem]{Lemma}

\newtheorem{definition}[theorem]{Definition}

\def\C{{{\mathbb C}}}

\def\T{{{\mathbb T}}}

\def\N{{{\mathbb N}}}

\def\sskip{\vskip 0.3 cm\noindent}

\newcommand{\QED}{\hfill$\Box$}

\newcommand{\Affines}[1]{\mathcal{A}(X)}
\newcommand{\Smooths}[1]{\mathcal{S}(X)}
\newcommand{\Rects}[1]{\mathcal{R}(X)}

\def\rd{{\rm d}}
\def\udisc{{\overline \Delta}}

\begin{document}

\maketitle

\begin{abstract}
We give a counterexample to a conjecture of S.E. Morris by showing that
there is a compact plane set $X$ such that $R(X)$ has no non-zero,
bounded point derivations but
such that $R(X)$ is not weakly amenable.
We also give an example of a separable uniform algebra $A$ such that every
maximal ideal of $A$ has a bounded approximate identity but such that
$A$ is not weakly amenable.
\end{abstract}

\section{Introduction}

The notion of weak amenability for commutative Banach algebras was originally
introduced
in \cite{BCD}. A commutative Banach algebra $A$ is said to be {\it weakly amenable} 
if there are no non-zero, continuous derivations from $A$ into any commutative Banach $A$-bimodule. 
However, as shown in \cite{BCD}, it is enough to check the dual module $A'$:  
$A$ is weakly amenable if and only if there are no non-zero, continuous
derivations from $A$ into $A'$. 

The condition of weak amenability sits between the stronger condition of
amenability (as defined by B.E. Johnson in
\cite{J}) and the weaker condition that the commutative Banach algebra $A$
have no non-zero, bounded point derivations (point derivations may be regarded 
as derivations into $1$-dimensional, commutative Banach $A$-bimodules).
There are examples in \cite{BCD} of commutative, semisimple 
Banach algebras showing that no two of these three conditions are
equivalent. (See also, for example, Sections 2.8 and 4.1 of \cite{D} 
for further discussion of this area and many illustrative examples.)

For uniform algebras, the differences between these three conditions
are less clear. Sheinberg's theorem (\cite[Theorem 5.6.2]{D}) 
tells us that
a uniform algebra $A$ on a compact, Hausdorff space $X$ is amenable if and only 
if $A=C(X)$ (the uniform algebra of all continuous, complex-valued functions
on $X$). In this case we say that the uniform algebra $A$ is
{\it trivial}.

No non-trivial, weakly amenable 
uniform algebras are known. There are examples of non-trivial, separable 
uniform algebras which have no non-zero, bounded point derivations (\cite{W})
and  even examples which have no non-zero point derivations at all (\cite{C}). 
In this paper we show that such uniform algebras need not be weakly amenable.

We will assume some familiarity with the standard relationships between 
bounded approximate identities, strong boundary points,
peak points and point derivations for uniform 
algebras (especially for $R(X)$).
Details may be found in the standard texts on uniform algebras
(for example
\cite{Browder}, \cite{Gamelin} or \cite{Stout}).
\sskip
{\bf Notation.} We denote the closed unit disc in $\C$ by
$\udisc$. More generally, $\Delta(a,r)$, ${\overline\Delta}(a,r)$
denote respectively the
open, closed, discs of radius $r$ centred on $a$.
For a bounded, complex-valued function $f$
defined on a non-empty set $S$ we
shall denote by $|f|_S$ the uniform norm of $f$ on $S$, that is
$$|f|_S = \sup\{|f(x)|: x\in S\}.$$
For a compact plane set $X$, $R_0(X)$ denotes the set of restrictions to
$X$ of rational functions with poles off $X$. Thus the standard uniform algebra
$R(X)$ is the uniform closure of $R_0(X)$ in $C(X)$.
\sskip
Suppose that $\mu$ is a measure on a compact plane set $X$ such that the bilinear functional 
defined on $R_0(X) \times R_0(X)$ by
$$(f,g)\mapsto\int_X f'(x)g(x) d \mu (x)$$
is bounded.
Then, as noted in \cite{F}, we may extend by continuity to $R(X) \times R(X)$
and obtain a continuous derivation $D$
from $R(X)$ to $R(X)'$ satisfying, for $f$ and $g$ in $R_0(X)$,
$$D(f)(g)=\int_X f'(x)g(x) d \mu (x).$$
Moreover, all continuous derivations from $R(X)$ to $R(X)'$ have this form.
Such a derivation is the zero derivation if and only if the measure $\mu$ annihilates
$R(X)$.

In \cite{Morris}, S.E. Morris investigated the extent to
which continuous derivations from $R(X)$ may be expressed as some form of 
integral of bounded point derivations. 
Morris conjectured 
that whenever $R(X)$ had no non-zero, bounded
point derivations then $R(X)$ must be weakly amenable.
We shall give a counterexample to this conjecture by constructing a 
suitable compact plane set
$X$ such that $R(X)$ has no non-zero, bounded point derivations but
such that $R(X)$ is not weakly amenable. (Our example will be a \lq Swiss
cheese'.)
In fact we shall prove the following theorem.

\begin{theorem}
For each $C>0$ there is a compact plane set $X$ obtained by deleting from the 
closed
unit disc a countable union of
open discs such that the unit circle $\T$ is a subset of $X$,
$R(X)$ has no non-zero, bounded point derivations, but for all 
$f$, $g$ in $R_0(X)$,
$$\left \vert \int_{\T} f'(z) g(z) ~\rd z \right \vert \leq C |f|_X |g|_X.$$
\end{theorem}

Let $X$ be one of the compact plane sets given by Theorem 1.1. 
In view of the discussion above
there is a non-zero, continuous derivation $D$
from $R(X)$ into $R(X)'$ satisfying, for all $f$ and $g$ in $R_0(X)$,
$$D(f)(g)=\int_{\T} f'(z) g(z) ~\rd z.$$
Thus $R(X)$ is not weakly amenable.
Note, however, that the map $f \mapsto f'|\T$ does not give a continuous derivation from
$R(X)$ into $L^1(\T)$,
as the polynomial functions $f(z)=z^n$ show. 

The problem (for compact plane sets $X$) 
of whether $R(X)$ can be weakly amenable without being $C(X)$
remains open. If you allow derivations into complete-metrizable modules then
the answer is negative: as shown in \cite{F}, $R(X)$ is 
\lq weakly-(F)-amenable' if and only if it is trivial.
\sskip
\noindent

\section{Construction of the Swiss cheeses}

To construct our Swiss cheeses we require several lemmas.

The first lemma is similar to some standard estimates
using Cauchy's theorem but has
an extra twist using Fubini's theorem.

\begin{lemma}
Let $D_n$ be a sequence of open discs in $\C$ (not necessarily
pairwise disjoint) whose closures are contained in the open unit disc.
Set
$X = \udisc\setminus \bigcup_{n=1}^{\infty}{D_n}$.
Let $s_n$ be the distance from $D_n$ to $\T$ and let
$r_n$ be the radius of $D_n$.
 Let $f$ and $g$ be in $R_0(X)$. 
Then

$$\left |{\int_{\T} f'(z) g(z) \rd z}\right | \leq
4\pi |f|_X |g|_X \sum_{n=1}^\infty
\frac{r_n}{s_n^2}.$$
\end{lemma}
\noindent
{\bf Proof.}
Let $\varepsilon>0$.
For $N \in \N$ we set
$$X_N={\overline \Delta}(0,(N+1)/N))\setminus \bigcup_{n=1}^N D_n.$$
Then we may choose $N$ large enough that $f$ and $g$ have no poles in $X_N$
and
$$4\pi \frac{N+1}{N}|f|_{X_N} |g|_{X_N}  \leq (4 \pi + \varepsilon) |f|_X |g|_X.\eqno{(1)}$$

Using the boundary of $X_N$, and ignoring any isolated points
(these may exist since the deleted discs can overlap), we obtain as usual (see, for example,  
pages 28-29 of \cite{BD})
a contour made up of
circular arcs which has winding number $1$ around all points of the interior of $X_N$
and $0$ around all points of $\C \setminus X_N$. We may split this contour into
$\gamma_1$ and $\gamma_2$ where $\gamma_1$ is the circle centered on the origin
and radius $(N+1)/N$, while $\gamma_2$ is the rest. Note that $\gamma_2$ then
satisfies
$$\int_{\gamma_2} \frac{|\rd \zeta|}{{\rm dist}(\zeta,\T)^2} \leq
\sum_{n=1}^N \frac{2\pi r_n}{s_n^2}.\eqno{(2)}$$

By Cauchy's formula, we have
$$f'(z)=h_1(z)+h_2(z)$$
where  (for $j=1,2$)
$$h_j(z)=\frac{1}{2\pi i} \int_{\gamma_j} \frac{f(w)}{(w-z)^2}~\rd w.$$
We now have
$$\int_{\T} f'(z) g(z) ~\rd z =
\int_{\T} h_1(z) g(z) ~\rd z +\int_{\T} h_2(z) g(z) ~\rd z.$$
We estimate the moduli of these two integrals separately.

The estimate for the $h_2$ integral is fairly standard:
by (2) we have, for $z \in \T$,
$$|h_2(z)| \leq |f|_{X_N} \sum_{n=1}^N \frac{r_n}{s_n^2}$$
and so
$$\left | {\int_{\T} h_2(z) g(z) ~\rd z} \right |
\leq 2\pi |f|_{X_N} |g|_{X_N} \sum_{n=1}^\infty \frac{r_n}{s_n^2}.\eqno{(3)}$$

The estimate for the $h_1$ integral requires an application of Fubini's theorem:
we have
$$\int_{\T} h_1(z)g(z) \rd z =
 \int_{\T} \left(\frac{1}{2\pi i}{\int_{\gamma_1} \frac{f(w)}{(w-z)^2} ~\rd w}\right) g(z) ~\rd z$$
$$=\frac{1}{2\pi i}
\int_{\gamma_1} f(w)\left({\int_{\T} \frac{g(z)}{(w-z)^2}~\rd z}\right) \rd w. \eqno{(4)}$$
However we may apply Cauchy's theorem to the function $g(z)/(w-z)^2$
to obtain
$$\int_{\T} \frac{g(z)}{(w-z)^2}~\rd z = - \int_{\gamma_2}  \frac{g(z)}{(w-z)^2}~\rd z$$
and so, by (2),
$$\left|{\int_{\T} g(z) /(w-z)^2~\rd z}\right| \leq |g|_{X_N}
 \sum_{n=1}^N 2\pi \frac{r_n}{s_n^2}.$$
Substituting this into equation (4) gives us
$$\left | {\int_{\T} h_1(z)g(z) \rd z}\right | \leq
2\pi \frac{N+1}{N} |f|_{X_N} |g|_{X_N} \sum_{n=1}^\infty \frac{r_n}{s_n^2}.\eqno{(5)}$$

Combining (1), (3) and (5) gives
$$\left |{\int_{\gamma} f'(z) g(z) \rd z}\right | \leq
(4\pi + \varepsilon) |f|_X |g|_X \sum_{n=1}^\infty
\frac{r_n}{s_n^2}.$$
Letting $\varepsilon$ tend to $0$ gives the result.
\QED

The next two results enable us to eliminate bounded point derivations by first working 
on compact sub-discs of $\Delta$.

\begin{lemma}
Let $X$ be a compact subset of $\udisc$.
Suppose that there is a sequence of real numbers $R_n \in (0,1)$ such that
$R_n \to 1$ and, for each $n$,
$R(X\cap{\overline\Delta}(0,R_n))$ has no non-zero, bounded point derivations. Then
$R(X)$ has no non-zero, bounded point derivations.
\end{lemma}
\noindent
{\bf Proof.}
It is clear that all points of $X \cap \T$ are peak points, and so there are no
non-zero point derivations at all at these points.
For the remaining points, the fact that $R(X)$ has no non-zero, bounded point derivations
follows from the local nature of this property, using (for example) Hallstrom's theorem 
(\cite[Theorem 1]{H}).
\QED

\sskip
The next lemma may be proved quickly by elementary arguments. It is also, however, an 
immediate consequence of Hallstrom's theorem.

\begin{lemma}
Let $Y$, $Z$ be compact plane sets with $Y \subseteq Z$.
If $R(Z)$ has no non-zero, bounded point derivations, then the same is true for $Y$.
\end{lemma}

Finally we quote a result proved by Wermer \cite{W} (see also \cite{H}).

\begin{proposition}
Let $D$ be a closed disc in $\C$ and let $\varepsilon>0$.
Then there is a sequence of open discs
$U_k \subseteq D$ such that 
$R(D\setminus \cup_{k=1}^\infty U_k)$ has no non-zero,
bounded point derivations but such that the sum of the radii of the discs $U_k$ is
less than $\varepsilon$.
\end{proposition}

We are now ready to construct the Swiss cheese we want.

\noindent
{\bf Proof of Theorem 1.1.}
Let $C>0$. Set $R_n=(n-1)/n$. By Proposition 2.4 we may 
choose open discs 
$U_{n,k} \subseteq {\overline \Delta}(0,R_n)$ 
such that $R({\overline \Delta}(0,R_n) \setminus
\cup_{k=1}^\infty{U_{n,k}}$ has no non-zero, bounded point derivations and
such that the sum (over $k$) of the radii of the discs $U_{n,k}$ is less than
$C(1-R_n)^2/(2^{n+3}\pi)$.

Set
$$X=\udisc \setminus \bigcup_{n,k} U_{n,k}.$$
Lemmas 2.2 and 2.3 show that $R(X)$ has no non-zero, bounded point derivations.
For the rest,
we may enumerate all the discs $U_{n,k}$ as $D_1,D_2,D_3,\dots$ and then apply
Lemma 2.1 to obtain the required estimate on the integral.
\QED

We conclude this section by discussing some additional properties of the 
uniform algebra $R(X)$ and the derivation $D: R(X) \rightarrow R(X)'$
constructed using Theorem 1.1 (for some $C>0$).
It is easily seen that this derivation is {\it cyclic}. 
This means that
it satisfies $(D(f))(g)= - (D(g))(f)$ for all $f$, $g$ in $R(X)$
(or, equivalently, $(D(f))(1)=0$ for all $f \in R(X)$).
\sskip
{\bf Notation.}
Let $A$ be a uniform algebra on a compact, Hausdorff space 
$X$ and let $E$ be a closed subset of $X$. 
We denote the ideal of functions in $A$ vanishing identically on $E$ by $I(E)$.

In this setting it is standard that 
the restriction algebra 
$A|E$ is isomorphic to $A/I(E)$ and so is a Banach function
algebra using the quotient norm. 

Taking $A$ to be our example of $R(X)$ above, the derivation $D$ constructed shows 
that $A$ is not weakly amenable. 
We now show that the Banach function algebra $A|\T$ is not weakly amenable
either. Since the derivation above is defined on elements of $R_0(X)$ using
only their values on $\T$ we would expect to be able to 
use this to define 
a derivation on $A|\T$. However, some caution is needed: so far we only know,
for $f$, $g$ in $R(X)$, that $|(D(f))(g)| \leq C |f|_X|g|_X$. For
$f$, $g$ in $R_0(X)$, though, we have 
$|D(f)(g)| \leq 2\pi |f'|_\T |g|_\T$.
Since, as we mentioned above, $D$ is cyclic we also have
$|D(f)(g)| \leq 2\pi |f|_\T |g'|_\T$.
It now follows by continuity, for $f$, $g$ in $R(X)$,
that if at least one of $f$, $g$ is in $I(\T)$ then $D(f)(g)=0$.
Thus, using some elementary functional analysis, 
we may define a non-zero, continuous
derivation $\tilde D$ from the Banach function algebra
$A|\T$ to $(A|\T)'$ by
$(\tilde D (f|\T))(g|\T) = (D(f))(g)$.

This shows that $A|\T$ is not weakly amenable. Since weak amenability passes
to quotients (by, for example, \cite[Proposition 2.1]{G}), 
it follows that any other uniform algebra with
a restriction isomorphic to $A|\T$ must also fail to be weakly amenable.
This observation will be crucial in the next section.

\section{Peak points without weak amenability}

The main result of this section will be the following theorem.

\begin{theorem}
There exists a uniform algebra $B$ 
whose maximal ideal space is metrizable (so $A$ is separable)
and such that every point of 
the maximal ideal space of $B$ is a peak point 
(equivalently, every maximal ideal
has a bounded approximate identity), but such that $B$ is not weakly 
amenable.
\end{theorem}

To construct this algebra, we will need a slightly modified version of the 
example constructed above, ensuring that there is a good supply of 
functions in the algebra vanishing identically on the circle. 
We will then use a version of Cole's systems of root extensions
(see \cite{C} and also \cite{F2})
to obtain the result.

We begin with some definitions 
and lemmas that we need. In particular we will prove
a result concerning the preservation 
of regularity properties by certain types of algebra extension,
including Cole's root extensions. 
(See \cite{F2} for some other results of this type.)

\sskip
{\bf Notation.} 
Let $A$ be a uniform algebra on a compact
space $X$.
For $x\in X$, we denote by $M_x$ and $J_x$ the ideals of
functions in $A$ vanishing at $x$, and in a neighbourhood of
$x$, respectively.
\sskip
We now recall the definition of a point of continuity (introduced
in \cite{FS}).

\begin{definition}
Let $A$ be a uniform algebra on a
compact Hausdorff space $X$ and let
$x\in X$. We say that $x$ is a {\sl point of continuity} for $A$
if there is no point $y$ of $X \setminus \{x\}$ satisfying $M_x\supseteq J_y$.

Equivalently, $x$ is a point of continuity for $A$ if and only if, for every 
compact set $E \subseteq X\setminus \{x\}$, there is a function $f$ in $A$ 
such that $f(x)=1$ and $f(E) \subseteq \{0\}.$ 

The algebra $A$ is {\it regular on $X$} if every point of 
$X$ is a point of continuity for $A$. In this case, if $X$ is the maximal ideal space of 
$A$ then $A$ is {\it regular}.
\end{definition}

See \cite{F3} for details of the connections between points of continuity and 
Jensen measures.

In the case where $E=\emptyset$, the next result follows from the work of 
Cole (\cite{C}), as every representing measure must then be a Jensen measure.
Recall also that for metrizable spaces strong boundary points and peak points 
coincide. In particular, this will be the case when we prove Theorem 3.1.

\begin{lemma}
Let $A$ be a uniform algebra on a compact space $X$ and let $x$ be a point of
continuity for $A$. Suppose that $E$ is a compact subset of $X\setminus\{x\}$
such that $\{f^2:f\in I(E)\}$ is dense in $I(E)$. Then $x$ is a strong
boundary point for $A$.
\end{lemma}
\noindent
{\bf Proof.} 
Note first, for convenience, that an easy application of the Mittag-Leffler inverse 
limit theorem 
(\cite[Corollary A.1.25]{D})
shows that $I(E)$ has a dense subset consisting of
functions with roots (in $A$, and hence in $I(E)$) of order $2^k$ for all $k\in \N$. 
Call this set of functions $S$.

Let $U$ be an open neighbourhood of $x$. We show that there
is a function $f \in A$ with $|f|_{X \setminus U} < 1/3$, $|f(x)|>2/3$
and $|f|_X < 2$. The result then follows from Gonchar's criterion 
(\cite[Corollary 7.20]{Stout}).

Since $x$ is a point of continuity for $A$ we may choose a function 
$g \in A$ such that $g(x)=1$, $g(X\setminus U) \subseteq \{0\}$
and $g(E) \subseteq 0$. Set $M=\|g\|_X$.

Choose $k \in \N$ such that $M^{2^{-k}} < 2$ and then choose $\delta \in (0,1)$
such that $\delta^{2^{-k}} < 1/3$, $(M+\delta)^{2^{-k}}<2$ and 
$(1-\delta)^{2^{-k}} > 2/3$.
Now choose $h \in S$ with $|h-g|_X < \delta$ and choose
$f \in A$ with $f^{2^k}=h$.
It is clear that the function $f$ has the required properties.
\QED
\sskip
The next Lemma shows that in some sense the property of being a point of continuity is
preserved by certain types of algebra extension. (Cole's systems of
root extensions are of this type.)
Note that we make no assumptions on the maximal ideal spaces of $A$ and $B$ 
here.

\begin{lemma}
Let $A$, $B$ be uniform algebras on, respectively, compact spaces $X$, $Y$
and suppose that there are a continuous surjection $\pi: Y \longrightarrow X$
and $x_0\in X$ satisfying:

(a) for all $f \in A$, $f\circ \pi$ is in $B$;

(b) the points of the fibre $F=\pi^{-1}(x_0)$ are separated by
a family of functions in $B$ each of which takes only finitely many 
different values on this fibre $F$ (so in particular $F$ is totally
disconnected).

If $x_0$ is a point of continuity for $A$ 
then every point of the fibre $F$ is a point of continuity for $B$.
\end{lemma}

{\bf Proof.}
Let $x \in F$ and $y \in Y\setminus \{x\}$.
We show that there exists $g$ in $B$ vanishing on some neighbourhood of
$y$ but with $g(x) \neq 0$.

Firstly, consider the case where $y$ is not in $F$. Then $\pi(y) \neq x_0$.
Since $x_0$ is a point of continuity for $A$ we may choose an $f \in A$ which vanishes on an open neighbourhood $U$ of $\pi(y)$ and with $f(x_0) \neq 0$.
Take $g=f\circ \pi$. Then $g$ vanishes on the open neighbourhood
$\pi^{-1}(U)$ of $y$ and $g(x) \neq 0$ as required.

It is now easy to see that for every compact set $E \subseteq Y \setminus F$
there is a function $h$ in $B$ vanishing on $E$ and with $h(x) \neq 0$.

We now consider the remaining case, where $y \in F$. By the assumption (b)
there is a function $b \in B$ taking only finitely many values on $F$
and with $b(x) \neq b(y)$.
Applying a polynomial to $b$ if necessary we may assume that 
$b$ takes only the values $0$ and $1$ on $F$, $b(y)=0$ and $b(x)=1$.

Choose an open neighbourhood $V$ of $F$ such that 
$b(V) \subseteq \Delta(0,1/3) \cup \Delta(1,1/3)$.
By Runge's theorem we may choose a sequence of polynomials
$p_n$ converging uniformly on $\Delta(0,1/3) \cup \Delta(1,1/3)$
to $\chi_{\Delta(0,1/3)}$.
Choose $h \in B$ vanishing on $Y \setminus V$ and with $h(x) \neq 0$.
Then the functions $h p_n(b)$ converge uniformly on $Y$ to a function
$g$. We then have that $g$ is in $B$, $g(x)=h(x) \neq 0$
and $g$ vanishes identically on the neighbourhood 
$b^{-1}(\Delta(0,1/3))$ of $y$, as required. \QED
\sskip
We are now ready to prove Theorem 3.1. In the interests of conciseness we 
leave some details to the reader. For further details on
systems of root extensions
the reader may wish to consult
\cite{C}, \cite[Section 19]{Stout} and, for a similar construction,
\cite[Section 5]{F2}.
\sskip
{\bf Proof of Theorem 3.1.}
Using similar methods to those of \cite{F3} (based on
K\"orner's version of McKissick's lemma, \cite{K}) we may easily modify 
the Swiss cheeses constructed in the proof of Theorem 1.1 to ensure, in 
addition, that
every point of $X\setminus \T$ is a point of 
continuity for $R(X)$. 

For such a compact plane set $X$ we set $A=R(X)$, and form a Cole
system of root extensions as follows. We first adjoin
square roots to a countable dense set of functions in $I(\T)$.
This ensures that a copy of 
$\T$ persists
in the maximal ideal space of the new uniform algebra.
We next adjoin square roots to a countable dense subset of
the ideal of functions in the new uniform algebra vanishing on this copy
of $\T$. Repeating this process we obtain a sequence of uniform algebras,
and we take the direct limit of this sequence.
This produces a uniform algebra $B$ on a compact, metrizable space $Y$ and a 
continuous surjection $\pi$ from $Y$ onto $X$
with the following four properties.
\begin{enumerate}
\item[(i)]
Conditions (a) and (b) of Lemma 3.4 hold for 
all $x_0$ in $X$.
\item[(ii)]
Set $E=\pi^{-1}(\T)$. Then $\pi|E$ is a homeomorphism
from $E$ onto $\T$.
\item[(iii)]
There is a norm-decreasing linear map $T$ from $B$ onto $A$
(obtained by averaging over fibres) such that
$T(f\circ \pi)=f$ for all $f \in A$.
\item[(iv)]
With $E$ as above,
$\{f^2:f \in I(E)\}$ is dense in $I(E)$.
\end{enumerate}

We claim that the uniform algebra $B$ has the desired properties.
First note that, for $g \in B$, we have
$g|E =(T(g)|\T)\circ \pi$. It follows easily that the restriction algebras
$A|\T$ and $B|E$ are isomorphic.
By the comments at the end of Section 2 
this shows that $B$ is not weakly amenable.

We wish to show that every point of $Y$ is a peak point for $B$.
For points in $E$ this is elementary as we may use $f\circ \pi$ for
appropriate peaking functions in $A$.
The remaining points lie in the fibres over points of continuity for
$A$ and so the result follows from Lemmas 3.3 and 3.4.

\QED

\section{Open questions}

We conclude with some open questions.
\begin{enumerate}
\item
Let $X$ be a compact plane set. Suppose that $R(X)$ is weakly amenable.
Must $R(X)=C(X)$?
\item
Are there any non-trivial, weakly amenable uniform algebras?
\item
Let $X$ be a compact plane set. 
Suppose that $R(X)$ has no non-zero, bounded point derivations. 
Can there ever be a measure $\mu$ on a compact subset $E$ of
$X$ such that the map $f \mapsto f'|E$ gives a non-zero, bounded 
derivation from
$R_0(X)$ into $L^1 (E,\mu)$?
\end{enumerate}

\sskip

{\sf  School of Mathematical Sciences

 University of Nottingham

 Nottingham NG7 2RD, England

 email: Joel.Feinstein@nottingham.ac.uk}
\sskip
2000 Mathematics Subject Classification: 46J10, 46H20
\end{document}